\documentclass[notitlepage,11pt]{amsart}
\usepackage{graphics}

\newtheorem{theorem}{Theorem}[section]

\theoremstyle{definition}
\newtheorem{definition}[theorem]{Definition}
\newtheorem{example}[theorem]{Example}

\theoremstyle{remark}

\numberwithin{equation}{section}

\DeclareSymbolFont{AMSb}{U}{msb}{m}{n}
\DeclareMathSymbol{\N}{\mathbin}{AMSb}{"4E}
\DeclareMathSymbol{\Z}{\mathbin}{AMSb}{"5A}
\DeclareMathSymbol{\R}{\mathbin}{AMSb}{"52}
\DeclareMathSymbol{\Q}{\mathbin}{AMSb}{"51}
\DeclareMathSymbol{\I}{\mathbin}{AMSb}{"49}
\DeclareMathSymbol{\C}{\mathbin}{AMSb}{"43}

\newcommand{\be}{\begin{equation}}
\newcommand{\ee}{\end{equation}}
\newcommand{\bea}{\begin{eqnarray*}}
\newcommand{\eea}{\end{eqnarray*}}

\begin{document}

\title{Bijective Quasi-Isometries of Amenable Groups}

\author{Tullia Dymarz}
\address{Department of Mathematics, University of Chicago, Chicago, Illinois 60637}
\email{tullia@math.uchicago.edu}
\thanks{The author was supported in part by NSERC PGS B}


\date{\today}



\begin{abstract}
Whyte showed that any quasi-isometry between non-amenable groups is a bounded distance from a bijection. In contrast this paper shows that for amenable groups, inclusion of a proper subgroup of finite index is never a bounded distance from a bijection.
\end{abstract}

\maketitle

\section{Introduction}

In his book on infinite groups, Gromov \cite[page 23]{Gr} asked whether inclusion of finite index subgroups
$F_m \longrightarrow F_n$ of two free groups is a bounded distance from a bi-Lipschitz map (i.e. a bijective quasi-isometry).
Papasoglu
answered this question affirmatively in \cite{Pa}.
A more general question asked in \cite[page 107]{H} is whether any two groups which are quasi-isometric always have a
bijective quasi-isometry between them? No counterexamples have been found. 

Whyte showed in \cite{W} that any quasi-isometry between \emph{non-amenable} groups is a bounded distance from a bijection. 
In contrast, we will show that for \emph{amenable} groups, inclusion of a finite index proper subgroup is never a bounded distance from a bijection (see Theorem \ref{tullia}). 
However, if the subgroup admits an ``$n$-to-$1$'' self quasi-isometry (where $n$ is the index of the subgroup) one can always compose this self map with the subgroup inclusion map
to get a quasi-isometry that is a bounded distance from a bijective quasi-isometry. So for such groups commensurability implies that
there does exist a bijective quasi-isometry.

\section{Preliminaries}

\begin{definition} A map between metric spaces $f:X\longrightarrow Y$ is a \emph{quasi-isometry} 
if there exist $C,K \geq 0$ such that for all $x,y \in X$  
\begin{eqnarray}\label{qi}
  -C+\frac{1}{K}d(x,y)\leq d(f(x),f(y)) \leq Kd(x,y) +C 
\end{eqnarray}
and there also exist $g:Y\longrightarrow X$ satisfying \ref{qi} such that $f\circ g$ and $g\circ f$ are a bounded 
distance from the identity. In this case $X$ and $Y$ are said to be \emph{quasi-isometric.}
\end{definition}

\begin{definition}A metric space $X$ is a \emph{uniformly discrete space of bounded geometry} (UDBG space)
if 
\begin{itemize}
\item there is an $\epsilon>0$ such that, for all $x,y\in X$, $d(x,y)<\epsilon \Rightarrow x=y$ (i.e. $X$ is uniformly discrete), and
\item  for any $r>0$ there is a bound $M_r$ on the size of any $r$-ball in $X$ (i.e. $X$ has bounded geometry).
\end{itemize}
\end{definition}

\begin{example} A finitely-generated group with the word metric $d(x,y)=||x^{-1}y||$ is a
main example of a UDGB-space. The norm is defined with respect to a specific generating set so two different generating sets may give metric spaces which are not isometric but they are quasi-isometric. 
\end{example}

\begin{example} Any proper metric space admitting a properly discontinuous isometric cocompact group action is 
also a UDBG-space.
\end{example}


To define the notion of an \emph{amenable} space we need some notation. Let $|S|$ denote the cardinality of the set $S$ and $\partial S=\{x\in X\ \mid \ 0<d(x,S)\leq 1\}$ denote the boundary of $S$.

\begin{definition}[Folner Criterion]
A UDBG-space is \emph{amenable} if there exists a sequence of finite subsets $\{S_i\}$ with 
\bea
\lim_{i\to \infty} \frac{|\partial S_i|}{|S_i|}=0
\eea
Such a sequence is called a \emph{Folner sequence}. 
\end{definition}

We say a finitely generated group amenable if it is amenable as a UDBG-space with the word metric. Since amenability is preserved by quasi-isometries this is a well defined notion.

It will also be useful to review some conventions of big ``O'' notation. 
\begin{definition}
We say 
\begin{itemize}
\item$f(i) = O(g(i))$ if there exist $C,K > 0$, such that $|f(i)|\leq C|g(i)|$ for all $i\geq K$.
\item $O(f(i)) < O(g(i))$ if $g(i) \neq O(f(i))$
\item  $O(g(i))=O(f(i))$ if both $g(i)= O(f(i))$ and $f(i)= O(g(i))$. 
\end{itemize}
\end{definition}
One important property we will use is that for $f,g \geq 0$
\bea O(f(i)+g(i))=O(f(i))\ \textrm{if} \ g(i)=O(f(i)).
\eea
Using this notation, the statement of amenability can be rephrased as follows: 
\\

\emph{ $X$ is amenable if there exists a sequence of finite sets $S_i\subset X$ such that $O(|\partial S_i|) < O(|S_i|)$.}

\section{Uniformly finite homology}

The uniformly finite homology groups $H_i^{uf}(X)$ were first introduced by Block-Weinberger in \cite{S}. Only $H_0^{uf}(X;\Z)$ is needed here, but
for a detailed discussion of uniformly finite homology, see \cite{W} or \cite{S}.  

For a UDBG-space $X$, let $C_0^{uf}(X)$ denote the vector space of infinite formal sums of the form
\bea
c=\sum_{x\in X}a_xx
\eea
 where there exists $M_c > 0$ such that $|a_x|\leq M_c$. 
Let $C_1^{uf}(X)$ denote the vector space of infinite formal sums of the form
\bea
c=\sum_{x,y\in X}a_{(x,y)}(x,y)
\eea
where there exist $M_c>0$ such that $|a_{(x,y)}|<M_c$ and $R_c>0$ such that $a_{(x,y)}=0$ if $d(x,y)>R_c$. 
Define a boundary map by
\bea
\partial : C_1^{uf}(X) &\longrightarrow & C_0^{uf}(X)\\
(x,y)& \longmapsto & y-x
\eea
and extending by linearity. Then we let 
\bea
H_0^{uf}(X)=C_0^{uf}(X)/\partial(C_1^{uf}(X))
\eea
 
Some important facts about $H_0^{uf}(X)$ we will not prove here (see \cite{W}) are
\begin{itemize}
\item  if $X$ and $Y$ are quasi-isometric then $H_0^{uf}(X)\cong H_0^{uf}(Y)$.
\item  if $X$ is infinite then $H_0^{uf}(X)$ is a vector space over $\R$. (When $X$ is finite  $H_0^{uf}(X)\cong \Z$)
\end{itemize}


\begin{definition} Any subset $S\subset X$ defines a class in $H_0^{uf}(X)$, denoted $[S]$, where [S] is the class of the chain $\sum_{x\in S}x$. We call $[X]$ the \emph{fundamental class} of $X$ in $H_0^{uf}(X)$.
\end{definition}

Using uniformly finite homology, Kevin Whyte developed in \cite{W} a test to determine 
when a quasi-isometry between UDBG spaces is a bounded distance from a
bijection.
\begin{theorem}[Whyte]\label{Whyte}\cite{W} Let $f:X \to Y$ be a quasi-isometry between
UDBG-spaces. Then there exists a bijective map 
a bounded distance from $f$ if and only if $f_{\ast}([X])=[Y]$. (Here $f_{\ast}([X])=[\sum_{x\in X}f(x)]$)

\end{theorem}
For non-amenable spaces we have the following theorem:
\begin{theorem}[Block-Weinberger]\label{Shmuel-Whyte}\cite{S}
Let $X$ be a UDBG-space. Then the following are equivalent:

\begin{itemize}
\item $X$ is non-amenable.

\item $H^{uf}_0(X)=0$ 

\item there exists $c= \sum_{x\in X} a_xx$ with $a_x > 0$ such that $[c]=0$ in $H^{uf}_0(X)$
\end{itemize}
 
\end{theorem}

Some of the motivation behind Whyte's Theorem \ref{Whyte} was that combined with Theorem \ref{Shmuel-Whyte} it implies:\\
\\
\emph{ Any quasi-isometry between finitely generated non-amenable 
groups is a bounded distance from
a bijection.}
\\
\\
 We can also use Theorem \ref{Whyte} to investigate 
quasi-isometries of amenable groups. To use Theorem \ref{Whyte} we need to be able to check when a chain in $c\in C_0^{uf}(X)$ represents the zero class in $H_0^{uf}(X)$. 
The following theorem gives such a criterion.

\begin{theorem}[Block-Weinberger]\label{need}\cite{S} (Theorem 7.6 in \cite{W})
Let $X$ be a UDBG-space, and let $c=\sum_{x\in X}a_xx \in C_0^{uf}(X)$. Then we have $[c]=0 \in H_0^{uf}(X)$ if and only if for any Folner sequence $\{S_i\}$, 
\bea
        \left| \sum_{x\in S_i}a_x \right|=O(|\partial S_i|).
\eea 
\end{theorem}
We now show how Whyte's criterion can be used to show that subgroup inclusion for amenable groups is not a
bounded distance from a bijection.
 
\begin{theorem}\label{tullia}
Let $G$ be an amenable group with proper subgroup $H$ of finite index, i.e. $[G:H]=n>1$. Then the
inclusion map $i:H \hookrightarrow  G$ is not a bounded distance from a bijective map.
\end{theorem}
\begin{proof}
Using Theorem \ref{need} we show that the chain $c=\sum_{x\in G\setminus H} x$
gives a nonzero class in $H_0^{uf}(G)$, that is $[c]=[G]-[H]\neq 0$. To this end let $\{S_i\}$ be any Folner sequence for $G$. 
Now $G=\bigcup_{k=1}^n g_kH$ and 
\bea
\sum_{k=1}^{n} |S_i\setminus g_kH|=(n-1)|S_i|
\eea
 so
\bea
   O(|S_i|)=O((n-1)|S_i|)=O(\sum_{k=1}^{n} |S_i\setminus g_kH|)=O(|S_i\setminus g_{k_i}H|)
\eea
for some $k_i$.
Let $F_i=g_{k_i}^{-1}S_i$. These sets also form a Folner sequence, since left multiplication by $g_{k_i}^{-1}$ is an 
isometry. Now $|g_{k_i}^{-1}S_i\setminus H|=|S_i\setminus g_{k_i}H|$. This gives us 
\bea 
 O(F_i\setminus H)= O(|F_i|) > O(|\partial F_i|).
\eea
So for the chain $c$, 
\bea
 \left|\sum_{x\in F_i}a_x \right|=|F_i\setminus H|=O(|F_i|) \neq O(|\partial F_i|),\\
\eea
and so $[c]\neq 0$.
\end{proof}

The following is a shorter proof of Theorem \ref{tullia} suggested by Weinberger.

\begin{proof}
  For the inclusion map $i$ to be a bijection, Theorem \ref{Whyte} tells us that we need 
\bea
[H]=i_{\ast}([H])=[G].
\eea
 But $[G]=n[H]$, so in fact we need $[H]=n[H]$. Now $G$ is amenable so 
$[H]\neq 0$. And since $H^{uf}_0(G)$ is torsion free $[H]=n[H]$ only if $n=1$.
\end{proof}

Theorem \ref{tullia} is actually a corollary of a more general result.

\begin{theorem}\label{second}
If $\phi:H\longrightarrow G$ is a homomorphism of amenable groups with finite index 
image, $[G:\phi(H)]=n$, and finite kernel, $ |\phi^{-1}(0)|=k $, then $\phi$ is 
a bounded distance from a bijection if and only if $n=k$.
\end{theorem} 
\begin{proof}
As above, to get a bijection we need $\phi_{\ast}([H])=[G]$. But 
$\phi_{\ast}([H])=k[\phi(H)]$ and $[G]=n[\phi(H)]$. Now $G$ is an amenable group, so
$[\phi(G)]\neq 0$, giving us $n=k$.
\end{proof}

\section{``$n$-to-$1$'' self Quasi-Isometries}
In this section we show how one can use subgroup inclusion to get bijective quasi-isometries between certain groups.
\begin{definition}
We call $f:X\longrightarrow X$ an \emph{``$n$-to-$1$'' self quasi-isometry}
if $f$ is a quasi-isometry and $|f^{-1}(x)|=n$. In this case 
$f_{\ast}([X])=n[X]$ . 
\end{definition}
Theorem \ref{second} suggests that we may be able to ``fix'' subgroup inclusion $i:H \rightarrow G$ of
index $n$ by precomposing with an $n$-to-$1$ quasi-isometry $f:G\rightarrow G$. Then $i\circ f:H \longrightarrow G$ is a new quasi-isometry that is a 
bounded distance from a bijection, since
\bea
(i\circ f)_{\ast}([H])=i_{\ast}(n[H])=n[H]=[G].
\eea
This leads to the question:
\\

\emph{Which amenable groups admit ``$n$-to-$1$'' quasi-isometries?}
\\
\\

\begin{example}
It is easy to define an ``$n$-to-$1$'' quasi-isometry  of $\Z$
\bea
\phi: & \Z\rightarrow\Z  \\
      &  k \mapsto  \lfloor \frac{k}{n}\rfloor
\eea
 where $\lfloor \frac{k}{n} \rfloor$ denotes the greatest
integer less than or equal to $\frac{k}{n}$. This is an ``$n$-to-$1$'' map of $\Z$ which is a $(\frac{1}{n},1)$ quasi-isometry. One can extend
this idea to get an ``$n$-to-$1$'' map on $\Z^m$ by applying the above map to one of the coordinates. 
\bea
\phi(k_1,k_2,\cdots ,k_m)=(\lfloor \frac{k_1}{n} \rfloor,k_2,\cdots,k_m)
\eea
\end{example}

We now consider another class of examples. 
\begin{example}
The \emph{solvable Baumslag Solitar groups} are given by the presentation 
\bea
BS(1,m)=\left<a,b\mid aba^{-1}=b^m\right>.
\eea

We can view $BS(1,m)$ as a union of cosets of the subgroup $\left< b \right>\cong \Z$. 
By identifying each coset with $\Z$ we can define an ``n-to-1'' map in a similar way as we do for $\Z$ 
\bea
f_{\alpha}:\alpha b^i \mapsto \alpha b^{ \lfloor \frac{i}{n}\rfloor}
\eea
where $\alpha$ is the coset representative. Picking a set $C$ of coset representatives gives us an ``n-to-1'' map of $BS(1,m)$
\bea
f_{C}:g \mapsto b^{ \lfloor \frac{i}{n}\rfloor}\alpha
\eea
where $g=\alpha b^i$ for some $\alpha \in C$.

 A priori, if we choose random coset representatives for each coset we may not get a quasi-isometry of $BS(1,m)$. It turns out that we can identify the cosets with $\Z$ in such a way so that the resulting map is actually a quasi-isometry.

In order to understand how to pick the coset representatives it is useful to review some ideas from \cite{FM}.
The group $BS(1,m)$ acts properly discontinuously and cocompactly on a
metric  $2$-complex $X_m$, which is a warped product of a tree $T_m$ and $\R$.
The tree $T_m$ is a $(m+1)$-valent directed tree with one ``incoming edge'' and $m$
``outgoing edges'' at each node. There is a natural projection $X_m \to
T_m$. The inverse image of a coherently oriented line, (a bi-infinite path in $T_m$ respecting the orientation) is a hyperbolic
plane. Any time we refer to these embedded hyperbolic planes we will identify them with the upper half plane model of $H^2$. The inverse image of a vertex is a horocycle, (called a branching 
horocycle). We can pick the basepoint of $X_m$ to lie on a branching horocycle.  

 $X_m$ is actually the universal cover of a complex $C_m$ (see \cite{FM} for details) whose fundamental group is $BS(1,m)$. The fundamental domain can be thought of as a ``horobrick'' $h_m \subset H^2$  defined by the region bounded by $0\leq x \leq n$ and $1\leq y \leq m$ so that the top of the horobrick has length $1$ and the bottom has length $m$ (see figure \ref{brick}).

\begin{figure}
  \begin{center}
    \includegraphics{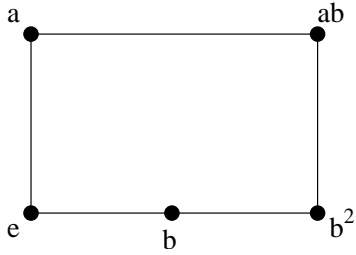}
  \end{center}
  \caption{the horobrick $h_2$ which is the fundamental domain corresponding to the element $a$ in $BS(1,2)$}
  \label{brick}
\end{figure}

We can define a quasi-isometry $i:BS(1,m) \rightarrow X_m$ by mapping e, the identiy of $BS(1,m)$, to the basepoint of $X_m$ and extending equivariantly. 
Then we can view $BS(1,m)$ as embedded in $X_m$ where group elements lie on branching horocycles and each branching horocycle contains a coset of $\left<b \right>$. Elements which differ only by the 
generator $a$ are distance $log(m)$ apart and lie on adjacent branching horocycles.
 Since $i$ is a quasi-isometry it has a coarse inverse. Let $j: X_m \rightarrow BS(1,m)$ be the coarse inverse of $i$ which maps each fundamental domain in $X_m$ to the unique element of $BS(1,m)$ in that domain. (i.e. each horobrick is mapped to the element in the upper left corner.)
Any map $f:BS(1,m) \rightarrow BS(1,m)$ gives us a map $i\circ f \circ j : X_m \rightarrow X_m$ which is a 
quasi-isometry of $X_m$ if and only if $f$ is a quasi-isometry of $BS(1,m)$.
 When convenient we will make no distinctions between $f$ and $i\circ f \circ j$.

The key idea is to ``line up'' all of the cosets so that our map $f$, when restricted to each hyperbolic plane in $X_m$, is a bounded distance from the quasi-isometry
\bea
 \phi: (x,y)\longmapsto (\frac{1}{n}x,y).
\eea

To this end  we need to consider another projection $X_m \rightarrow H^2$. From \cite{FM} we know that there exists a unique map $\rho_m: X_m \rightarrow H^2$ with the following properties:
\begin{itemize}
\item $\rho_m$ takes horocycles to horocycles
\item $\rho_m$ is an isometry when restricted to each hyperbolic plane in $X_m$
\item  $\rho_m$ is normalized to take the base point of $X_m$ to the point $(x,y)=(0,1)$
\end{itemize}
Let $l=\{ (0,y) \in H^2 \}$ be the $y$-axis and consider $T=\rho^{-1}(l)$. $T$ intersects each branching horocycle at exactly one point. This will be our reference point which we will call ``0''. 
Since each branching horocycle ``contains'' a coset of $\left<b\right>$, pick the coset representative for this coset to be a group element $\alpha$ which lies closest to ``0''. 
(There may be two such elements in some cases). Note that for each $\alpha$ the distance between $\alpha$ and $T$ is at most one. Our map $i\circ f_\alpha \circ j$ is bounded distance from the map $\phi$ when restricted to each branching horocycle. Because all of the $\alpha$'s are a uniformly bounded distance from $T$ we have that for each hyperbolic plane $Q$ in $X_m$ the total map $f_C$ restricted to $Q$ is a bounded distance from $\phi$. So by the rubberband principle $i\circ f_C \circ j$ is a quasi-isometry of $X_m$ and so $f_C$ is an ``n-to-1'' quasi-isometry of $BS(1,m)$.

\end{example}

We now consider one criterion for when a group does admit an ``$n$-to-$1$'' quasi-isometry. If $G$ contains a subset $G'$ such that 
\bea
G=\bigsqcup_{i=1}^n G'g_i.
\eea
and there exists a bijective quasi-isometry 
\bea
f:G'\longrightarrow G
\eea
then $f$ extends to an  ``$n$-to-$1$'' self quasi-isometry of $G$ given by
\bea
f':G &\longrightarrow &G\\
   g'g_i &\longmapsto & f(g') \ \ (g'\in G')
\eea
This holds in particular if $G'$ is a subgroup of $G$ with $|G:G'|=n$ and $f:G' \rightarrow G$ is an isomorphism.

 If it were possible to find an ``$n$-to-$1$'' self quasi-isometry
for all amenable groups, then we would have a bijective quasi-isometry between any two commensurable, amenable groups.

\specialsection*{Acknowledgments}

I'd like to thank Kevin Whyte for useful conversations, Shmuel Weinberger, and especially my advisor Benson Farb for many helpful suggestions.

\bibliographystyle{amsalpha}

\end{document}